\documentclass[a4paper,12pt]{article}
\usepackage{amsmath,amsfonts,amssymb,bbm,units,amsthm} 

\usepackage[latin2]{inputenc}
\usepackage[T1]{fontenc}

\usepackage{geometry}
\newtheorem{theorem}{Theorem}[section]

\newtheorem{proposition}{Proposition}[section]
\newtheorem{corollary}{Corollary}[section]

\newtheorem{remark}{Remark}[section]

\newcommand{\ud }{{\ \rm d}}
\newcommand{\rf}[1]{\eqref{#1}}

\newcommand{\bbfR}{{\mathbb R}}

\newcommand{\vf}{{\varphi}}
\newcommand{\ve}{{\varepsilon }}
\newcommand{\se}{{\text{\rm{e}}}}
\def\qed{\hfill{$\Box$}}

\title{Asymptotic stability of singular solution to nonlinear heat equation}
\author{Dominika Pilarczyk\\
\rule{0pt}{15pt}\\
\small Instytut Matematyczny, Uniwersytet Wroc\l awski,\\
\small pl. Grunwaldzki 2/4, 50-384 Wroc\l aw, Poland \\
\small e-mail: \small \texttt {dpilarcz{@}math.uni.wroc.pl}}

\date{\rule{0pt}{15pt}\\ \today}

\begin{document}

\maketitle


\begin{abstract}
{\footnotesize
In this paper, we discuss the asymptotic stability of singular steady states of the nonlinear heat equation $u_t=\Delta u+u^p$
in weighted  $L^r$-- norms.}
\end{abstract}

\section{Introduction}
\setcounter{equation}{0}

\label{ce}
It is well-known that the behaviour for large $t$ of solutions of the Cauchy problem
\begin{align}
      \label{h.e.}&u_t=\Delta u+u^p,\\
      \label{d.}  &u(x,0)=u_0(x)
\end{align}
depends on the value of the exponent $p$ of the nonlinearity. Let us first recall the critical value of $p=p_F=1+\nicefrac{2}{n}$ called the Fujita exponent which borders the case of a finite-time blow-up for all positive solutions (for  $p\leqslant p_F$) and the case of the existence of some global bounded positive solutions (if $p>p_F$). It is also known that the Sobolev exponent $p_S=\frac{n+2}{n-2}$
 is critical for the existence of positive steady states that is classical solutions $\psi \in C_0( {\bbfR ^n})$
of the elliptic equation $\Delta \psi +\psi ^p=0\quad {\rm on} \quad {\bbfR ^n} .$
Such solutions exist only if $p \geqslant p _S $ (see e.g. \cite{Chen}, \cite{Gidas} ). Moreover, for $p \geqslant p _S$ there is a one parameter family of radial positive steady states $\psi _k,\ k>0$, given by 
\begin{equation}\label{psi k1}
      \psi _k(x)=k\psi_1(k^{\frac{p-1}{2}}|x|),
\end{equation}          
where $\psi _1$ is the unique radial stationary solution with $\psi _1(0)=1$, which is stricly decreasing in $|x|$ and satisfies $\psi _1(|x|)\rightarrow 0$ as $|x|\rightarrow \infty $ (see \cite{QS}).

Another important exponent 
\begin{equation*}    
     p _{JL}=\frac{n-2\sqrt{n-1}}{n-4-2\sqrt{n-1}} \quad {\rm for }\quad n\geqslant 11,
\end{equation*}     
appeared for the first time in \cite{JL} where the authors studied problems with the non\- linearities of the form $f(u)=a(1+bu)^p$ for some $a, b>0$. It is also connected with a change in stability property of  positive steady states defined in \rf{psi k1}. Indeed, Gui {\it et al.} \cite{GNW} proved that for $p <p_{JL}$, all positive stationary solutions $\psi _k$ are unstable in any re\-asonable sense, while for $p\geqslant p_{JL}$ they are "weakly asymptotically stable" in a weighted $L^\infty-$norm. Results on the asymptotic stability of zero solution to \rf{h.e.}-\rf{d.} can be found in \cite{Q08} and in the references given there.

Let us recall that for $p \geqslant p_{JL}$ the family of the positive equilibria $\psi _k$, $k>0$, forms a simply ordered curve. Furthermore, this curve connects the trivial solution if $k \rightarrow 0$ and the singular steady state for $k \rightarrow \infty $, which exists for $p>p_{st}=\nicefrac{n}{(n-2)}$ in dimensions $n \geqslant 3$ 
and has the form $v_{\infty }(x)=L |x|^{\nicefrac{-2}{(p-1)}}$ with a suitably chosen constant $L$ (see \rf{v inf} below). It is also known (\cite{QS}) that if $p_S\leqslant p <p_{JL}$ the graphs of the steady states $\psi _k$, $0<k<\infty $, intersect the graph of $v_{\infty }$, whereas for $p \geqslant  p _{JL}$ we have $\psi _k<v_{\infty }$, $(0<k<\infty ).$ 

Our main goal in this note is to prove asymptotic stability of the singular stationary solution $v_\infty $ in suitable weighted $L^r-$spaces using estimates of a fundamental solution to a parabolic equation with singular coefficients \cite{LS,MS}.

\section{Results and comments}
\setcounter{equation}{0}

It can be directly checked that for $p>p_{st}=\nicefrac{n}{(n-2)}$ and $n\geqslant 3$ equation \rf{h.e.} has the singular stationary solution of the form
\begin{equation}\label{v inf}
       v_{\infty }(x)=L |x|^{-\frac{2}{p-1}}=\Bigg(\frac{2}{p-1}\Big(n-2-\frac{2}{p-1}\Big)\Bigg)^\frac{1}{p-1} |x|^{-\frac{2}{p-1}},
\end{equation}
which plays the central role in this paper.

In particular, problem \rf{h.e.}--\rf{d.} with a nonnegative initial datum $u_0$, which is bounded and below singular steady state $v_\infty $, has the global in time classical solution (see \cite[Th. 20.5 (i)]{QS} and \cite[Th. 1.1]{PY}). Moreover, following Galaktionov \& Vazquez \cite[Th. 10.4 (ii)]{Gal-Vaz}, we may generalize that result and prove that if  $0\leqslant u_0(x)\leqslant v_\infty (x)$ and $u_0(x)\not \equiv v_\infty (x)$, then the limit function $u(x,t)=\lim_{N\rightarrow \infty }u^N(x,t),$
where $u^N=u^N(x,t)$ is the solution of the problem
\begin{align*}
   &u_t=\Delta u+u^p, \quad u(x,0)=\min \{u_0 (x),N\},
\end{align*}   
solves \rf{h.e.} and $u(\cdot ,t)\in L^\infty (\bbfR^n)$ for all $t>0$. By those reasons, in the theorems below we always assume that $u$ is the nonnegative solution to the initial value problem \rf{h.e.}--\rf{d.} with the initial datum $u_0$ satisfying 
\begin{equation}\label{a.1}
     0\leqslant u_0(x) \leqslant v_\infty (x).
\end{equation}

In order to show the asymptotic stability of the steady state $v_\infty $ we linearize \rf{h.e.} around $v_\infty $. Denoting by $u=u(x,t)$ the nonnegative solution to \rf{h.e.}--\rf{d.} and introducing $w=v_\infty -u$, we obtain
\begin{equation}\label{l.e. w}              
             w_t=\Delta w+\frac{\lambda }{|x|^2}w-\big[ (v_\infty -w)^p-v_\infty ^p+pv_\infty ^{p-1}w \big], 
\end{equation}
where 
\begin{equation}\label{lambda}
         \lambda =\lambda (n,p)= \frac{2p}{p-1}\Big(n-2-\frac{2}{p-1}\Big).
\end{equation}             
Next, we use estimates of the fundamental solution of the linear heat equation with singular potential
\begin{equation}\label{lin}
         u_t=\Delta u +\frac{\lambda }{|x|^2}u, \quad x\in \bbfR^n , \quad t>0
\end{equation} 
obtained recently by Liskevich \& Sobol \cite{LS}, Milman \& Semenov \cite{MS} (see also Moschini \& Tesei \cite{MT}). As the consequence of the Hardy inequality, it is crucial in that reasoning to assume that $\lambda \leqslant \frac{(n-2)^2}{4}$ in equation \rf{lin}. Coming back to the perturbed equation \rf{l.e. w} and using the explicit form of $\lambda (n,p)$ in \rf{lambda}, we obtain by direct calculation (see Remark \ref{1: ex lin e} for more details) that the inequality $\lambda (n,p)\leqslant \frac{(n-4)^2}{4}$ is valid if 
\begin{equation}\label{a.2}
      p\geqslant p_{JL}=\frac{n-2\sqrt{n-1}}{n-4-2\sqrt{n-1}} \quad \textrm{for }\quad n\geqslant 11.
\end{equation}
By this reason, we limit ourselves to the exponent $p$ of the nonlinearity in \rf{h.e.} satisfying \rf{a.2}. The exponents mentioned above are ordered as follows: $p_F < p_{st} <p_S < p_{JL}$.

We introduce the parameter $\sigma $ which plays a crucial role in our reasoning by the formula
\begin{equation}\label{a.3}
       \sigma =\sigma(n,p)= \frac{n-2}{2}-\sqrt{\frac{(n-2)^2}{4}-\frac{2p}{p-1}\bigg(n-2-\frac{2}{p-1}\bigg)}.
\end{equation}
It is worth pointing out that $\sigma (n,p) >\nicefrac{2}{(p-1)}$ if $p>p_{st}$ and $n>2$. Moreover, the number $\sigma (n,p)$ has the property $2\sigma (n,p)<n$. Let us also notice that $\sigma (n,p)$ appears in a hidden way in the papers of Pol\'a\v cik, Yanagida, Fila, Winkler (see {\it e.g.}\cite{PY1}, \cite{FWY}), because it is the sum of the constant $\nicefrac{2}{(p-1)}$ and $\lambda_1$, where $\lambda _1$ is one of the root of the quadratic polynomial $z^2-(n-2-2L)z+2(n-2-L)$, given explicitly by the formula  
\begin{equation*}
\lambda_1=\frac{1}{2}\bigg(n-2-2L-\sqrt{(n-2-2L)^2-8(n-2-L)} \bigg),
\end{equation*}
where $L$ is defined in \rf{v inf}.

Now we are in a position to formulate our first result on the convergence of the solutions towards the singular steady state.
\begin{theorem}\label{half l}
   Assume \rf{a.1}, \rf{a.2}, \rf{a.3}. Suppose, moreover, that there exist constants $b>0$ and $\ell \in \big(\sigma ,n-\sigma \big)$ such that 
\begin{equation*}
     v_\infty (x)-b|x|^{-\ell} \leqslant u_0(x)
\end{equation*}
  for all $|x|\geqslant 1$. Then
\begin{align}
     &\label{rel i} \sup_{|x|\leqslant \sqrt{t} }|x|^\sigma \big( v_\infty (x) - u(x,t) \big) \leqslant Ct^{-\frac{\ell -\sigma }{2}} \\
     \intertext{and}
     &\label{rel ii} \sup_{|x|\geqslant \sqrt{t}} \big( v_\infty (x) -u(x,t) \big) \leqslant Ct^{-\frac{\ell }{2}}.
\end{align}
for a constant $C>0$ and all $t\geqslant 1$.
\end{theorem}

Pol\'a\v cik \& Yanagida \cite[Th. 6.1]{PY} showed  that under the assumptions of Theorem \ref{half l} the pointwise convergence holds true, namely, $\lim_{t \rightarrow \infty}u(x,t)=v_\infty (x)$ for every $x \in \bbfR^n \setminus \{0\}$. More recently, Fila \& Winkler \cite{FW} proved the uniform convergence of solutions $u=u(x,t)$ toward a singular steady state on $\bbfR^n \setminus B_\nu (0)$, where $B_\nu (0)$ is the ball in $\bbfR^n$ with the center at the origin and radius $\nu $. Theorem \ref{half l}completes those results by providing optimal weighted decay estimates in the whole $\bbfR^n$.

\begin{remark}
{\rm
   Note that our calculations in the proof of Theorem \ref{half l} are valid for any $\ell \in (\nicefrac{2}{(p-1)}, n-\sigma )$, but for $\ell \in (\nicefrac{2}{(p-1)}, \sigma ]$ the right-hand side of inequality \rf{rel i} does not decay in time. 
}\qed
\end{remark}

We can improve Theorem \ref{half l} for $\ell=\sigma $ as follows.
\begin{theorem}\label{mth 2}
    Assume that \rf{a.1}, \rf{a.2} and \rf{a.3} are satisfied. Suppose that there exists a constant $b>0$ such that
\begin{equation*}
     v_\infty (x)-b|x|^{-\sigma  }\leqslant u_0 (x).
\end{equation*}
    Let, moreover,
\begin{equation*}
      \lim_{|x|\rightarrow \infty }|x|^{\sigma }\big( v_\infty (x)-u_0(x) \big) = 0.
\end{equation*}    
   Then 
\begin{align*}
     &\lim_{t\rightarrow \infty }\sup_{|x|\leqslant \sqrt{t}}|x|^\sigma \big(v_\infty (x)-u(x,t)\big)=0\\
     \intertext{\it and }
    &\lim_{t\rightarrow \infty }t^{\frac{\sigma }{2}}\sup_{|x|\geqslant \sqrt{t}}\big(v_\infty (x)-u(x,t)\big) =0.
\end{align*}
\end{theorem}

\begin{corollary}\label{small b} 
   Under the assumptions of Theorem \ref{half l} and Theorem \ref{mth 2}, respectively if, moreover,  $b$ is sufficiently small, we obtain
\begin{align}
     \label{sup u}   &\|u(\cdot ,t)\|_\infty \geqslant C t^{\frac{\ell -\sigma}{\sigma (p-1)-2}} \quad \textit{if } \quad \ell \in (\sigma , n-\sigma )\\
\intertext{\it for a constant $C>0$ and all $t\geqslant 1$ and }
     \label{sup u 2}&\lim_{t\rightarrow \infty }\|u(\cdot ,t)\|_\infty =+\infty  \quad \textit{if }\quad \ell=\sigma .
\end{align}
\end{corollary}

\begin{remark}
{\rm
Estimates from below of $\|u(\cdot ,t)\|_\infty $, similar to that stated in \rf{sup u}, were obtained by Fila {\it et al.} in \cite[Theorem 1.1.]{FWY1}, \cite[Theorem 1.1.]{FWY} and improved in  \cite[Theorem 1.1.]{FKWY} using matched asymptotics expansions. In Corollary \ref{small b}, we emphasize that this inequality is an immediate consequence of Theorem \ref{half l}.
}\qed
\end{remark}

\begin{remark}
{\rm
For $p>p_{JL}$ estimates \rf{rel i} and \rf{rel ii} seem to be optimal, because they imply the optimal lower bound \rf{sup u}, see \cite{FKWY}. On the other hand, for $p=p_{JL}$ the authors of \cite{FKWY1} obtained the logarithmic factor on the right-hand side of \rf{sup u}, which we are not able to see by our method.
}\qed
\end{remark}

Our next goal is to prove the asymptotic stability of $v_\infty $ in the Lebesgue space $L^2 (\bbfR^n )$.

\begin{theorem}\label{stab.2}
    Assume that \rf{a.1}, \rf{a.2} and \rf{a.3} are valid. 
\begin{itemize}
    \item[i)]  Suppose that $v_\infty -u_0 \in L^1(\bbfR^n )$ and $|\cdot |^{-\sigma } (v_\infty -u_0)\in L^1(\bbfR^n)$. Then 
\begin{equation}\label{1:  L1}
    \| v_\infty (\cdot ) - u(\cdot , t) \|_2 \leqslant Ct^{-\frac{n}{4}}\| v_\infty  -u_0  \|_1 + Ct^{-\frac{n-2\sigma }{4}}\| |\cdot |^{-\sigma }(v_\infty -u_0 )\|_1 .
\end{equation}    
    \item[ii)] Suppose that $v_\infty -u_0 \in L^2 (\bbfR^n )$. Then
\begin{equation*}
    \lim_{t \rightarrow \infty } \|v_\infty (\cdot )-u(\cdot ,t) \|_2 =0.
\end{equation*}
\end{itemize}
\end{theorem}

 Note, that $v_\infty \in L^2_{loc}(\bbfR^n )$ for every $p>p_F$. Here, this property of the singular solution $v_\infty $ is satisfied, because $p_{JL} >p_F$. 

Using the fact that the steady states $\psi_k$ defined in \rf{psi k1} are below the singular stationary solution $v_\infty $ for $p> p_{JL}$, we may rephrase Theorem \ref{stab.2} as follows.

\begin{corollary}\label{c. stab.2}
     Assume that \rf{a.1}, \rf{a.2} and \rf{a.3} are valid. Let $\psi_k$ be the stationary solutions \rf{psi k1} for some $k>0$. 
\begin{itemize}
    \item[i)]  Suppose that $\psi_k -u_0 \in L^1(\bbfR^n )$ and $|\cdot |^{-\sigma } (\psi_k -u_0)\in L^1(\bbfR^n)$. Then 
\begin{equation}\label{1: c L1}
    \| \psi_k (\cdot ) - u(\cdot , t) \|_2 \leqslant Ct^{-\frac{n}{4}}\| \psi_k  -u_0  \|_1 + Ct^{-\frac{n-2\sigma }{4}}\big\| \, |\cdot |^{-\sigma }(\psi_k -u_0 )\big\|_1 .
\end{equation}    
    \item[ii)] Suppose that $\psi_k -u_0 \in L^2 (\bbfR^n )$. Then
\begin{equation}\label{c lim L2}
    \lim_{t \rightarrow \infty } \| \psi_k (\cdot )-u(\cdot ,t) \|_2 =0.
\end{equation}
\end{itemize}
\end{corollary}

\begin{remark}
{\rm
Observe that Theorem \ref{stab.2} and Corollary \ref{c. stab.2} complete the results by Pol\'a\v cik \& Yanagida, who proved in \cite[Proposition 3.5]{PY} the stability estimate
\begin{equation*}
     \| \psi_k (\cdot ) -u(\cdot ,t)\|_2 \leqslant \|\psi_k -u_0\|_2 .
\end{equation*}
}\qed
\end{remark}

\section{Linear equation with a singular potential}
\setcounter{equation}{0}

In this section we recall the estimate from above of the fundamental solution of the equation $u_t=\Delta u+\lambda |x|^{-2}u$ obtained by Liskevich \& Sobol in \cite{LS} and by Milman \& Semenov \cite{MS}. Following those arguments, we define the weights $\vf_\sigma (x,t) \in C(\bbfR^n \setminus \{0\})$ as 
\begin{equation}\label{weights}
         \vf_\sigma (x,t)=
         \begin{cases}
            \big(\frac{\sqrt{t}}{|x|}\big) ^\sigma & {\rm if}\ |x| \leqslant \sqrt{t},\\
            1 & { \rm if}\  |x|\geqslant \sqrt{t} .
         \end{cases}  
\end{equation}

 \begin{theorem}{\rm \cite{LS,MS}}\label{kernel th}
       Let $Hu=\Delta u+\lambda |x|^{-2}u$. Assume that $0 \leqslant \lambda \leqslant \nicefrac{(n-2)^2}{4}$.
      The semigroup $\se^{-tH}$ of the linear operators generated by $H$ can be written as the integral operator with a kernel $\se^{-tH}(x,y)$, namely
\begin{equation*}
         \se^{-tH}u_0 (x)=\int _{\bbfR ^n}\se^{-tH}(x,y)u_0(y) \ud y.
\end{equation*}    
      Moreover, there exist positive constants $C>0$ and $c > 1$, such that for all $t>0$ and all $x,y \in \bbfR^n \setminus  \{0\}$
\begin{equation}\label{kernel}
         {0 \leqslant \se^{-tH}(x,y) \leqslant C\varphi_\sigma (x,t)\ \varphi_\sigma (y,t)\ G(x-y, c t)},
\end{equation}
      where $\sigma =\frac{n-2}{2}-\sqrt{\frac{(n-2)^2}{4}-\lambda }$, the functions $\vf_\sigma $ are defined in \rf{weights} (see also Remark \ref{R Phi} below) and $G(x,t)=\big(4\pi t\big)^{\nicefrac{-n}{2}}\exp( \nicefrac{-|x|^2}{4\pi t})$ is the heat kernel.
\end{theorem}

\begin{remark}\label{R Phi}
{\rm
   In fact, Milman \& Semenov in \cite{MS} used the more regular weight functions $\Phi _\sigma \in C^2(\bbfR^n \setminus \{0\})$, namely
\begin{equation*}
    \Phi_\sigma (x,t)=
         \begin{cases}
            \big(\frac{\sqrt{t}}{|x|}\big) ^\sigma & {\rm if}\ |x| \leqslant \sqrt{t},\\
            \frac{1}{2} & { \rm if}\  |x|\geqslant 2\sqrt{t} 
         \end{cases}  
\end{equation*}
   and $\frac{1}{2}\leqslant \Phi_\sigma (x,t) \leqslant 1$ for $\sqrt{t}\leqslant |x|\leqslant 2\sqrt{t}$. It can be checked directly that there exist positive constants $c$ and $C$ for which the inequalities 
\begin{equation*}
   c\vf_\sigma (x,t)\leqslant \Phi_\sigma (x,t) \leqslant C\vf_\sigma (x,t)
\end{equation*}
hold true, where $\vf_\sigma $ are defined by \rf{weights}. By this reason we are allowed to use the weights $\vf_\sigma $ instead of $\Phi_\sigma $.
}\qed
\end{remark}

The following theorem is the consequence of the estimates stated in \rf{kernel}. 

\begin{theorem}\label{w half l}
   Let the assumptions of Theorem \ref{kernel th} be valid. Assume that $p>1+\frac{2}{n-\sigma }$.
   Suppose that there exist $b>0$ and $\ell \in (\frac{2}{p-1},n-\sigma )$ such that a nonnegative function $w_0$ satisfies
\begin{alignat*}{2}
  &w_0(x)\leqslant b|x|^{-\frac{2}{p-1}} &\quad &\textit{for} \quad |x|\leqslant 1 ,\\
  &w_0(x)\leqslant b|x|^{-\ell} &\quad &\textit{for}\quad  |x|\geqslant 1.
\end{alignat*}
  Then
\begin{equation}\label{sup etH}
   \sup_{x\in \bbfR^n} \vf_\sigma ^{-1}(x,t) |\se^{-tH}w_0(x)|\leqslant Ct^{-\frac{\ell}{2}}
\end{equation}
  for a constant $C>0$ and all $t\geqslant 1$.
\end{theorem}

\begin{proof}
  First, for every fixed $x \in \bbfR^n $, we apply the estimate of the kernel $\se^{-tH}$ from Theorem \ref{kernel th} in the following way
\begin{equation*}
    \vf_\sigma ^{-1}(x,t) \big|\se^{-tH}w_0(x)\big| \leqslant C\int_{\bbfR^n} G(x-y,ct)\vf_\sigma (y,t)w_0(y) \ud y.
\end{equation*}
  Next, we split the integral on the right-hand side into three parts $I_1(x,t)$, $I_2(x,t)$ and $I_3(x,t)$ according to the definition of the weights $\vf_\sigma $ and the assumptions on the function $w_0$. Let us begin with $I_1(x,t)$:
\begin{equation*}
     \begin{split}
      I_1(x,t)&\equiv  C\int_{|y|\leqslant 1} G(x-y,ct)\vf_\sigma (y,t)w_0(y) \ud y \\
                 &\leqslant Cb t^{\frac{\sigma }{2}}\int_{|y|\leqslant 1} G(x-y, ct)|y|^{-\sigma -\frac{2}{p-1}} \ud y 
      \leqslant Cbt^{-\frac{n-\sigma }{2}},
      \end{split}
\end{equation*} 
because $G(x-y,ct)$ is bounded by $C t^{-\frac{n}{2}}$ and the function $|y|^{-\sigma -\frac{2}{p-1}}$ is integrable for $|y|\leqslant 1$ if $p>1+\nicefrac{2}{(n-\sigma )}$. 

We use the same argument to deal with 
\begin{equation*}
   \begin{split}
    I_2(x,t)&\equiv C\int_{1\leqslant |y|\leqslant \sqrt{t}}G(x-y,ct)\vf_\sigma (y,t)w_0(y) \ud y\\
     &\leqslant Cbt^{\frac{\sigma }{2}} \int_{1\leqslant |y|\leqslant \sqrt{t}} G(x-y,ct) |y|^{-\sigma -\ell} \ud y  
        \leqslant Cbt^{\frac{\sigma -n}{2}} \int_{1\leqslant |y|\leqslant \sqrt{t}}|y|^{-\sigma -\ell} \ud y \\
     &\leqslant Cbt^{-\frac{\ell}{2}} + Cbt^{-\frac{n-\sigma }{2}}.
   \end{split}
\end{equation*}

Finally, we estimate 
\begin{equation*}
    \begin{split}
     I_3(x,t)&\equiv C\int_{|y|\geqslant \sqrt{t}} G(x-y, ct)\vf_\sigma(y,t) w_0(y) \ud y \\
     &\leqslant Cb\int_{|y|\geqslant \sqrt{t}} G(x-y, ct)|y|^{-\ell} \ud y \leqslant Cbt^{-\frac{\ell}{2}},
     \end{split}
\end{equation*}
using the inequality $1\leqslant \big(\frac{\sqrt{t}}{|y|}\big)^\ell $ for $|y|\geqslant \sqrt{t}$ and the identity $\int_{\bbfR^n} G(x-y,ct) \ud y=1$ for $t>0$, $x\in \bbfR^n $.
Since $\ell \in (\nicefrac{2}{(p-1)},n-\sigma )$, we complete the proof of  \rf{sup etH}.
\end{proof}

\begin{theorem}\label{limit e -tH th}
    Assume that  $|\cdot |^\sigma w_0 \in L^\infty (\bbfR^n )$ and 
\begin{equation*}  
   \lim_{|x|\rightarrow \infty }|x|^\sigma w_0 (x)=0.
\end{equation*}   
   Then 
\begin{equation}\label{limit e -tH}
   \lim_{t\rightarrow \infty }t^{\frac{\sigma }{2}}\sup _{x\in \bbfR^n}\vf^{-1} _\sigma (x,t)|\se^{-tH}w_0(x)|=0.
\end{equation}
\end{theorem}

\begin{proof}
For every fixed $x\in \bbfR^n $ we use the estimate from Theorem \ref{kernel th} as follows
\begin{equation*}
   \vf_\sigma ^{-1}(x,t) \big| \se^{-tH}w_0(x)\big| \leqslant C\int_{\bbfR^n} G(x-y,ct)\vf_\sigma (y,t)w_0(y) \ud y.
\end{equation*}
We decompose the integral on the right-hand side according to the definition of $\vf_\sigma $ and we estimate each term separable. Substituting $y=z\sqrt{t}$ we obtain
\begin{equation*}
    \begin{split}
     I_1(x,t)&\equiv C\int_{|y|\leqslant \sqrt{t}} G(x-y,ct) \bigg( \frac{\sqrt{t}}{|y|} \bigg)^{\sigma }w_0(y) \ud y\\
     &=Ct^{-\frac{\sigma }{2}}\int_{|z|\leqslant 1 }G\bigg(\frac{x}{\sqrt{t}}-z,c\bigg)|z|^{-2\sigma} |\sqrt{t}z|^{\sigma }w_0(\sqrt{t}z) \ud z .
     \end{split}
\end{equation*}
Hence, 
\begin{equation*}
    t^{\frac{\sigma }{2}}\sup_{x\in \bbfR^n }I_1(x,t) \rightarrow 0 \quad \textrm{as} \quad t\rightarrow \infty
\end{equation*}
by the Lebesgue dominated convergence theorem, because $G\big(\frac{x}{\sqrt{t}}-z,c\big)$ is bounded and the function $|z|^{-2\sigma }$ is integrable for $|z|\leqslant 1$. By the assumption imposed on $w_0$, given $\ve >0$ we may choose $t$ so large that
\begin{equation*}
   \sup_{|y|\geqslant \sqrt{t} }|y|^\sigma w_0(y)<\ve .
\end{equation*}

Now, using the inequality $1\leqslant \big( \frac{\sqrt{t}}{|y|} \big)^\sigma $ for $|y|\geqslant \sqrt{t}$, we obtain 
\begin{equation*}
    \begin{split} 
     I_2(x,t)&\equiv \int_{|y|\geqslant \sqrt{t}} G(x-y,ct) w_0(y) \ud y 
     \leqslant t^{-\frac{\sigma }{2}}\int_{|y|\geqslant \sqrt{t}} G(x-y,ct)|y|^{\sigma }w_0(y) \ud y \\
     &\leqslant \ve t^{-\frac{\sigma }{2}}\int_{|y|\geqslant \sqrt{t}}G(x-y, ct) \ud y .
     \end{split}
\end{equation*}
Since $\int_{\bbfR^n } G(x-y,ct )\ud y=1$ for all $t>0$, $x\in \bbfR^n $ and since $\ve >0$ is arbitrary, we get 
\begin{equation*}
     t^{\frac{\sigma }{2}} \sup_{x\in \bbfR^n} I_2(x,t) \rightarrow 0 \quad \textrm{as} \quad t\rightarrow \infty .
\end{equation*}
\end{proof}

Let us defined the weighted $L^q$-norm as follows
\begin{equation*}
      \|f\|_{q,\vf_\sigma (t)}=\bigg( \int_{\bbfR^n} |f(x)\vf_\sigma ^{-1} (x,t)|^q \vf_\sigma ^2 (x,t) \ud x \bigg) ^{\frac{1}{q}} \quad \textrm{ for every}    \quad 1\leqslant q <\infty ,
\end{equation*}      
     and     
\begin{equation*}
    \|f\|_{\infty, \vf_\sigma (t)}=\sup_{x\in \bbfR^n} \vf_\sigma ^{-1}(x,t)|f(x)| \quad {\rm for }\quad q=\infty .
\end{equation*}
Note, that in particular for $q=2$, the norm $\| \cdot \|_{2, \vf_\sigma (t)} $ agrees with the usual $L^2$-norm on $\bbfR^n $. 

\begin{proposition}
Suppose that $1\leqslant q\leqslant \infty $. Then the following inequality holds true
\begin{equation}\label{w norm}
    \|\se^{-tH}w_0\|_{q,\vf _\sigma (t)}\leqslant Ct^{-\frac{n}{2}(\frac{1}{r}-\frac{1}{q})} \|w_0\|_{r,\vf_\sigma (t)}
\end{equation}
for every $1\leqslant r\leqslant q\leqslant \infty$ and all $t>0$.
\end{proposition}

\begin{proof}
The proof of estimate \rf{w norm} can be directly deduced from the reasoning by Milman \& Semenov. Indeed, in \cite[page 381]{MS}, we can find the inequality
\begin{equation*}
    \| \vf_\sigma ^{-1} \se^{-tH} \vf_\sigma f\|_{L^2\big(\bbfR^n , \vf_\sigma ^2 (x,t)  \ud x \big)} \leqslant Ct^{-\frac{n}{4}} \| f\|_{L^1\big(\bbfR^n , \vf_\sigma ^2 (x,t)  \ud x \big)}.
\end{equation*}     
Hence, substituting $\vf_\sigma f=w_0$ and using the definitions of the norm $\| \cdot \|_{q,\vf_\sigma (t)}$, we obtain \rf{w norm} with $q=2$ and $r=1$. This inequality together with 
\begin{equation*} 
     \|\se^{-tH}w_0\|_{1,\vf _\sigma (t) }\leqslant C\|w_0\|_{1,\vf_\sigma (t)},
\end{equation*}
stated in \cite[page 391]{MS}, imply \rf{w norm} for $q=1$ and every $1\leqslant r\leqslant 2$ by Riesz-Thorin interpolation theorem.
Moreover, the operator $\se^{-tH}$ is self-adjoint, so by duality the inequality 
\begin{equation*}
     \|\se^{-tH}w_0\|_{\infty , \vf_\sigma (t)}\leqslant Ct^{-\frac{n}{4}}\|w_0\|_{2, \vf_\sigma (t)}
\end{equation*}
holds true. The semigroup property $\se^{-tH}=\se^{-\nicefrac{t}{2}H}\se^{-\nicefrac{t}{2}H}$ leads to \rf{w norm} with $q=\infty $ and $r=1$.  Applying duality and Riesz-Thorin interpolation theorem once more, we complete the proof of \rf{w norm}. Let us emphasize at the end of this reasoning, that the inequalities \rf{w norm} are used by Milman \& Semenov in \cite{MS} to derive the kernel estimate \rf{kernel}.
\end{proof}
\section{Linearization around a singular steady state}
\setcounter{equation}{0}

Let $u$ be a solution of \rf{h.e.} with initial datum satisfying \rf{a.1}. We substitute 
\begin{equation*}      
    w(x,t)=v_\infty (x)-u(x,t)
\end{equation*}    
to get
\begin{equation}\label{h.e.1}   
      w_t=\Delta w+\frac{\lambda }{|x|^2}w-\big[ (v_\infty -w)^p-v_\infty ^p+pv_\infty ^{p-1}w \big],
\end{equation}
where $\lambda =\lambda (n,p)=\frac{2p}{p-1}(n-2-\frac{2}{p-1})$. Let us note that the last term on the right-hand side of equation \rf{h.e.1} is non positive, namely
\begin{equation*} 
(v_\infty -w)^p-v_\infty ^p \geqslant -pv_\infty ^{p-1}w,
\end{equation*}
which is the direct consequence of the convexity of the function $f(s)=s^p$. Indeed, since the graph of the function $f$ lies above all of its tangents, we have $f(s-h)-f(s)\geqslant -f'(s)h$ for all $s$ and $h$ in $\bbfR$.

The proofs of our results are based on the following elementary observation. If $w$ is a nonnegative solution of equation \rf{h.e.1} with the initial condition $w_0(x)\geqslant 0$, then 
\begin{equation*}
    0\leqslant w(x,t)\leqslant \se^{-tH}w_0(x)
\end{equation*}
with $Hw=\Delta w +\lambda (n,p)|x|^{-2}\, w $. Consequently, using the condition $0\leqslant u_0(x) \leqslant v_\infty (x)$ and the just-mentioned comparison principle we can write 
\begin{align}
    \label{v inf and e H 0}  & 0\leqslant v_\infty (x) -u(x,t)\leqslant \se^{-tH}\big(v_\infty(x)-u_0(x) \big)\\
    \intertext{or , equivalently,    }
    \label{v inf and e H}    &v_\infty (x) -\se^{-tH}\big( v_\infty (x)-u_0(x) \big) \leqslant u(x,t)\leqslant v_\infty (x).
\end{align}

\begin{remark}\label{1: ex lin e}
{\rm
    If  $n\geqslant 11$  and either $p \geqslant p_{JL}$ or $\frac{n}{n-2} <p <\frac{n+2\sqrt{n-1}}{n-4+2\sqrt{n-1}}$, then the linearized problem
\begin{equation*}
       \begin{split}  
       & w_t=\Delta w+\frac{\lambda }{|x|^2}w,\\
       &w(x,0)=w_0(x),
      \end{split}    
\end{equation*} 
    with $\lambda =\lambda (n,p) =\frac{2p}{p-1}(n-2-\frac{2}{p-1})$ has the unique solution.
Indeed, in the view of Theorem \ref{kernel th}, it is sufficient to show that
\begin{equation*}         
    \lambda (n,p)=\frac{2p}{p-1}\bigg(n-2-\frac{2}{p-1}\bigg) \leqslant \frac{(n-2)^2}{4}.
\end{equation*}    
Substituting $y=\nicefrac{1}{(p-1)}$, after elementary calculations, we arrive at the inequality
\begin{equation*}
    16y^2+(32-8n)y+n^2-12n+20 \geqslant 0
\end{equation*}    
which has the solution $y\in \big(-\infty , \frac{n-4-2\sqrt{n-1}}{4}\Big]\cup \Big[\frac{n-4+\sqrt{n-1}}{4}, +\infty \big)$. Moreover, if \ $n \geqslant 11$, then\  $\frac{n-4-2\sqrt{n-1}}{4}>0$ and if $n\in (2, 10)$, then $\frac{n-4-2\sqrt{n-1}}{4}<0$ and $\frac{n-4+2\sqrt{n-1}}{4}>0$.
These observations give us that $p \geqslant p_{JL}$ or $\frac{n}{n-2}<p\leqslant \frac{n+2\sqrt{n-1}}{n-4+2\sqrt{n-1}}$.
}\qed
\end{remark}

\section{Asymptotic stability of steady states}
\setcounter{equation}{0}

\begin{proof}[ Proof of Theorem \ref{half l}]
It suffices to use inequality \rf{v inf and e H 0} and to estimate its right-hand side by Theorem \ref{w half l}.
\end{proof}

\bigskip

\begin{proof}[Proof of Theorem \ref{mth 2}]
As in the proof of Theorem \ref{half l}, it is sufficient to use \rf{v inf and e H 0} together with Theorem \ref{limit e -tH th} substituting $w_0(x)=v_\infty (x)- u_0(x)$.
\end{proof}

\bigskip

\begin{proof}[Proof of Corollary \ref{small b}]
Since we have inequality \rf{v inf and e H}, it suffices to prove that
\begin{equation*}
    \sup_{x\in \bbfR^n}\big[v_\infty (x)-\se^{-tH}w_0(x)\big]\geqslant C(b) t^{\frac{\ell-\sigma }{\sigma (p-1)-2}}
\end{equation*} 
for $w_0=v_\infty -u_0$. Hence, inequality \rf{sup etH} from Theorem \ref{w half l} enables us to write
\begin{equation*}
    v_\infty (x)-\se^{-tH}\big(v_\infty (x)- u_0(x) \big)\geqslant v_\infty (x)-Cb \vf_\sigma (x,t)t^{-\frac{\ell}{2}}
\end{equation*}
for all $x\in \bbfR^n \setminus \{0\}$ and $t>0$.
Next, using the explicit form of the weights $\vf_\sigma $, we define the function
\begin{equation*}
    F(|x|,t)=v_\infty (|x|)-Cb\vf_\sigma (x,t) t^{-\frac{\ell }{2}}=
        \begin{cases}
         L|x|^{-\frac{2}{p-1}}-Cbt^{\frac{\sigma -\ell}{2}}|x|^{-\sigma } & {\rm for }\  |x| \leqslant \sqrt{t},\\
         L|x|^{-\frac{2}{p-1}}-Cbt^{-\frac{\ell }{2}} & {\rm for }\  |x|\geqslant \sqrt{t}.
      \end{cases}
\end{equation*}
  
An easy computation shows that the function $F$ has its maximum at 
\begin{equation*}          
    |x|=C(b ) t^{ \frac{\sigma - \ell}{2}\frac{p-1}{\sigma (p-1)-2}}
\end{equation*}    
and it is equal to
\begin{equation*}
   \max_{x\in \bbfR^n }F(|x|,t)=C(b) t^{\frac{\ell-\sigma }{\sigma (p-1)-2}}
\end{equation*}   
for some constant $C(b)\geqslant 0$. Hence, we get \rf{sup u}.

To obtain \rf{sup u 2}, we use the result from Theorem \ref{limit e -tH th}. It follows from \rf{limit e -tH} that for every $\ve >0$ there exists $T>0$ such that 
\begin{equation*}
     \big| \se ^{-tH}w_0 (x) \big| <\ve \vf_\sigma  (x,t) t^{-\frac{\sigma }{2}}
\end{equation*}
for all $x \in \bbfR^n \setminus \{0\}$ and $t>T$.
Hence, by \rf{v inf and e H}, we have
\begin{equation*}
      v_\infty (x)-\se^{-tH}\big(v_\infty (x)- u_0(x) \big)\geqslant v_\infty (x)-C \ve \vf_\sigma (x,t)t^{-\frac{\sigma }{2}}.
\end{equation*} 
Now, once more  using the explicit form of the weights $\vf_\sigma $, we consider the function 
\begin{equation*}
    G(|x|,t)=v_\infty (|x|)-Cb\vf_\sigma (x,t) t^{-\frac{\sigma }{2}}
    = \begin{cases}
         L|x|^{-\frac{2}{p-1}}-\ve |x|^{-\sigma } & {\rm for }\  |x| \leqslant \sqrt{t},\\
         L|x|^{-\frac{2}{p-1}}-\ve t^{-\frac{\sigma }{2}} & {\rm for }\  |x|\geqslant \sqrt{t}.
       \end{cases}
\end{equation*}
An elementary computations give us that the function $G$ attains its maximum at
\begin{equation*}          
    |x|=c \ve^{ -\frac{p-1}{\sigma (p-1)-2}}
\end{equation*}    
and 
\begin{equation*}
   \max_{x\in \bbfR^n }G(|x|,t)=C\ve^{-\frac{2}{\sigma (p-1)-2}}
\end{equation*}   
for some constant $C\geqslant 0$. Since $\sigma>\nicefrac{2}{(p-1)}$, we see that the maximum of the function $G$ diverges to infinity if $\ve $ tends to zero. This completes the proof of \rf{sup u 2}.
\end{proof}        

\bigskip

\begin{proof}[Proof of Theorem \ref{stab.2}(i)]
According to \rf{v inf and e H 0} it is enough to estimate the $L^2-$ norm of the expression $\se^{-tH}w_0$ for every $w_0$ satisfying two conditions: $w_0\in L^1 (\bbfR^n )$ and $| \cdot |^{-\sigma }w_0 \in L^1 (\bbfR^n )$. Applying \rf{w norm}, with $q=2$, $r=1$ and using the definition of the functions $\vf_\sigma (x,t)$, we may write
\begin{equation*}
     \begin{split}
     \|\se^{-tH} w_0\|_2 &\leqslant Ct^{-\frac{n}{4}}\|w_0\|_{1, \vf_\sigma (t)} =Ct^{-\frac{n-2\sigma }{4}}\int_{|x|\leqslant \sqrt{t}}w_0(x)|x|^{-\sigma }\ud x\\
                                   &+Ct^{-\frac{n}{4}}\int_{|x|\geqslant \sqrt{t}}w_0(x)\ud x  \leqslant Ct^{-\frac{n-2\sigma }{4}}\| w_0 |\cdot |^{-\sigma }\|_1 + Ct^{-\frac{n}{4}}\| w_0\|_1.
     \end{split}
\end{equation*} 
This establishes formula \rf{1: L1}.    
\end{proof}

\bigskip

\begin{proof}[Proof of Theorem \ref{stab.2}(ii)]
Again, by \rf{v inf and e H 0},  we only need to show that  
\begin{equation*}
     \lim_{t\rightarrow \infty } \|\se^{-tH}w_0\|_2 =0
\end{equation*}
for each $w_0 \in L^2 (\bbfR^n)$. Hence, for every $\ve >0$ we choose 
$\psi \in C_c^\infty ( \bbfR^n )$ such that $\|w_0 -\psi\|_2<\ve $. 
Using first the triangle inequality and next \rf{w norm}, with $q=2$ and $r=2$, we obtain
\begin{equation*}
     \begin{split}
     \| \se^{-tH}w_0 \|_2 &\leqslant \|\se^{-tH}(w_0 -\psi )\|_2 + \| \se^{-tH} \psi \|_2 \\
     &\leqslant C\ve + \| \se^{-tH} \psi \|_2.
     \end{split}
\end{equation*} 
Since the second term on the right-hand side convergence to zero as $t\rightarrow \infty $ by the first part of Theorem \ref{stab.2}, we get
\begin{equation*}
     \limsup_{t\rightarrow \infty }\| \se^{-tH}w_0\|_2 \leqslant C\ve.
\end{equation*}     
This completes the proof of Theorem \ref {stab.2} (ii), because $\ve >0$ can be arbitrary small.  
\end{proof}

\bigskip

\begin{proof}[Proof of Corollary \ref{c. stab.2}]
We linearize equation \rf{h.e.} around the positive steady state $\psi_k$ substituting $v=\psi_k -u$ to get
\begin{equation}\label{eq. v}
     v_t=\Delta v+p\psi_k^{p-1} v-\big( (\psi_k-v)^p - \psi_k^p + p\psi_k^{p-1}v \big).
\end{equation}
Once more, using the convexity of the function $f(s)=s^p$, let us notice that the expression $(\psi_k -v )^p-\psi_k^{p-1}v$ is nonnegative. Furthermore, $\psi_k <v_\infty $ as long as $p>p_{JL}$ and $n\geqslant 11$. Hence applying first, the comparison principle to the approximate problem
\begin{equation*}
     \begin{split}
      &v_t=\Delta +p \min \{ N^{p-1}, v_\infty ^{p-1} \}  v,\\
      &v(x,0)=v_0(x)
      \end{split}
\end{equation*}      
with the nonnegative initial datum and next, passing to the limit when $N$ tends to infinity, we get
\begin{equation*}
     0\leqslant \psi_k(x)-u(x,t)\leqslant \se^{-tH}\big( \psi_k(x) - u_0(x)\big) .
\end{equation*}
Now, \rf{1: c L1} and \rf{c lim L2} are the straightforward consequences of the reasoning used in the proof of Theorem \ref{stab.2}.     
\end{proof}

\section*{Acknowledgment}

The author wishes to express her gratitude to the anonymous referee for several helpful comments and to Jacek Zienkiewicz for many stimulating conversations and . This work is a part of the author PhD dissertation written under supervision of Grzegorz Karch.  


\end{document}